\documentclass[12pt, reqno]{amsart}
 
\usepackage{latexsym}

\newcommand{\boxt}{\Box\kern -6.3pt\raise .55pt
    \hbox{$\scriptstyle{\times}$}}

\newcommand{\del}{\partial}

\newcommand{\mapright}[1]{\smash{\mathop
   {\longrightarrow}\limits^{#1}}}

\newtheorem{thm}{Theorem}[section]
\newtheorem{lem}[thm]{Lemma}

\newtheorem{prop}[thm]{Proposition}
\newtheorem{cor}[thm]{Corollary}

\theoremstyle{definition}

\theoremstyle{remark}
\newtheorem{remark}[thm]{Remark}
\newtheorem{remarks}[thm]{Remarks}

\numberwithin{equation}{section}

\newcommand{\inv}{^{-1}} 
\newcommand{\ovl}{\overline}

\newcommand{\half}{\frac{1}{2}} 
\newcommand{\thalf}{\textstyle\frac{1}{2}}

                                     \newcommand{\g}{{\mathfrak g}}

\newcommand{\fs}{{\mathfrak s}}

\newcommand{\fso}{\mathfrak{so}} 
\newcommand{\fsp}{\mathfrak{sp}}
\newcommand{\fsl}{\mathfrak{sl}}

\newcommand{\og}{\ovl{g}}  \newcommand{\of}{\ovl{f}}

\newcommand{\C}{{\mathbb C}}

\newcommand{\om}{\omega}  
\newcommand{\al}{\alpha}
\newcommand{\be}{\beta}
\newcommand{\ga}{\gamma}
\newcommand{\ep}{\varepsilon}

\newcommand{\la}{\lambda} \newcommand{\La}{\Lambda}

\newcommand{\sig}{\sigma}
 
\newcommand{\vth}{\vartheta}
\newcommand{\ze}{\zeta}

                                            \newcommand{\A}{{\mathcal A}}
 
                                         \newcommand{\D}{{\mathcal D}}

\newcommand{\cH}{{\mathcal H}}
                                             \newcommand{\I}{{\mathcal I}}

\newcommand{\cK}{{\mathcal K}}

                                          \newcommand{\OO}{{\mathcal O}}

\newcommand{\cS}{{\mathcal S}} 
  
                                           \newcommand{\U}[1][]{{\mathcal U}_{#1}}

\newcommand{\Cl}{\mathrm{Cl}}

\newcommand{\gr}{\mathop{\mathrm{gr}}\nolimits}

\newcommand{\End}{\mathop{\mathrm{End}}\nolimits}
\newcommand{\Hom}{\mathop{\mathrm{Hom}}\nolimits}
\newcommand{\Ker}{\mathop{\mathrm{Ker}}\nolimits}

\newcommand{\Omin}{{\OO}_{\mathrm{min}}}

\newcommand{\gog}{\g\oplus\g}  
  
\newcommand{\Ug}[1][]{{\mathcal U}_{#1}({\mathfrak g})} 
 
\newcommand{\Sg}[1][]{S^{#1}({\mathfrak g})}  
 
\newcommand{\Kill}[2]{{\langle{#1},{#2}\rangle}_\g}

\newcommand{\bh}{{\mathbf h}}

\newcommand{\bq}{{\mathbf q}} \newcommand{\bqt}{{\mathbf q}_t}

\newcommand{\bs}{{\mathbf s}}

\newcommand{\RR}{{\mathcal R}}
\newcommand{\UgJ}{{\U/J}}
\newcommand{\Jla}{J^{\la}}   \newcommand{\Jh}{J^{\scriptscriptstyle\half}}
\newcommand{\Dla}{{\mathfrak D}^{\la}}   
            \newcommand{\Dhalf}{{\mathfrak D}^{\scriptscriptstyle\half}}

\newcommand{\starsub}[1][]{\star_{\scriptscriptstyle #1}}
\newcommand{\Pmux}[1][\cdot]{\{\mu_x,#1\}}
\newcommand{\Oreg}{\Omin^{*}}
\newcommand{\pep}{p_\ep}  \newcommand{\pmep}{p_{-\ep}}

\newcommand{\CP}{\mathbb{CP}}
\newcommand{\RP}{\mathbb{RP}}

\newcommand{\RPn}{\RP^{n}} \newcommand{\CPn}{\CP^n}
\newcommand{\gs}{\g^{\sharp}}

\newcommand{\ns}[1]{|\hskip -1.5pt |{#1}|\hskip -1.5pt |^2}

\begin{document}   

\title{Non-Local Equivariant Star Product on the Minimal Nilpotent Orbit}

\author{Alexander Astashkevich}
\address{Renaissance  Technologies} 
\email{ast@rentec.com}

\author{Ranee Brylinski}  
\thanks{Research supported in part  by NSF  Grant No.   DMS-9505055}   
\address{Department of Mathematics,
        Penn State University, University Park 16802}
\email{rkb@math.psu.edu}
\urladdr{www.math.psu.edu/rkb}

\dedicatory{Dedicated to Dmitri Fuchs on his 60th birthday}
\date{}

\begin{abstract} 
We construct a unique $G$-equivariant graded star product on  the algebra
$\Sg/I$ of polynomial functions on the minimal nilpotent  coadjoint
orbit $\Omin$  of $G$  where $G$ is a complex simple Lie group and 
$\g\neq\fsl(2,\C)$.  This   strengthens the result of Arnal, Benamor   and  Cahen.

Our main result is to compute, for $G$ classical, 
the star product  of a momentum function $\mu_x$ with any   function $f$. 
We find $\mu_x\star f=\mu_xf+\half\{\mu_x,f\}t+\La^x(f)t^2$. For
$\g$   different from $\fsp(2n,\C)$,  
$\La^x$ is not a differential operator. Instead $\La^x$ is the left quotient of an
explicit order  $4$ algebraic differential operator  $D^x$ by an order $2$ invertible
diagonalizable operator. Precisely,
$\La^x=-\frac{1}{4}\frac{1}{E'(E'+1)}D^x$ where $E'$ is a positive shift of the
Euler vector field. Thus $\mu_x\star f$ is not local in $f$.

Using $\star$ we construct a  
positive definite  hermitian inner product on $\Sg/I$. The Hilbert space completion of $\Sg/I$ is
then a unitary representation of $G$. This quantizes $\Omin$ in the sense of geometric
quantization and the orbit method.
\end{abstract}

\maketitle

\section{Introduction} 
\label{sec:intro}      
   
The fundamental problem in equivariant quantization is the $G$-equivariant
quantization of the coadjoint orbits of $G$, where $G$ is a simply-connected Lie
group.  In deformation quantization, there is a nice set of axioms for the star
product $\star$ and then  $G$-equivariance of $\star$ is a relation involving the
momentum functions $\mu_x$, $x\in\g$, where $\g=Lie(G)$. In fact, this 
amounts to $G$-equivariance of the
corresponding quantization map (see \S\ref{sec:2}).  

It was already recognized by Fronsdal (\cite{Fron}) that the locality axiom for 
star products must be modified in order to accommodate equivariance.
The locality axiom means, in either  the smooth or algebraic setting,
that the operators which define the star product  are bidifferential.

One could simply exclude any constructions that are  not  local. But this would  
cast aside  equivariant constructions  (such as  \cite[\S9, page
124]{Fron},\cite{CGsphere},\cite{CGreg},  and, as we show, \cite{ABC})
which are unique and very natural;   moreover   these retain a strong  flavor of
locality.   Figuring out what this ``flavor" is and how to
axiomatize it is a   very interesting problem. 
It seems to involve  ``pseudo-differential" operators. 

In this paper, we investigate the unique $G$-equivariant graded star product 
on the algebra $\RR$  associated to the   minimal (non-zero)  nilpotent  
coadjoint orbit $\Omin$ in $\g^*$, where  $\g$ is a simple complex Lie algebra
different from $\fsl(2,\C)$.
Here $\RR=\Sg/I$ is the algebra of polynomial functions on $\Omin$. 
The star product was constructed for  $\g$  different from $\fsl(n,\C)$ by Arnal,
Benamor,  and  Cahen in \cite{ABC}. We strengthen their result in our Proposition
\ref{prop:e+u},  after some preliminary work in \S\ref{sec:2}.
We find  an  analog of the Joseph ideal for $\g=\fsl(n,\C)$, $n\ge 3$. 
(There is a $1$-dimensional family of
candidates,  but   only one of them produces a star product with parity.)
We prove uniqueness whenever $\g\neq\fsl(2,\C)$.

To start off, we show (Proposition \ref{prop:4})  that the   star product of a
momentum function $\mu_x$ with any   function $f$ is the  three term sum
\begin{equation}\label{eq:1} 
\mu_x\star f=\mu_xf+\half\{\mu_x,f\}t+\La^x(f)t^2
\end{equation}
where $\La^x$ are   graded operators on $\RR$ of degree $-1$.
We compute $\La^x$ for $\g$ classical. For $\g=\fsp(2n,\C)$, we find 
(\S\ref{sec:5}) some familiar order $2$ differential operators (which appear in   
the Fock space model of the oscillator representation).  

Our main result (Theorem \ref{thm:Lax=})
is a formula for $\La^x$ when $\g$ is classical but different from
$\fsp(2n,\C)$, i.e., when  $\g=\fsl(n+1,\C)$ ($n\geq 2$) or  
$\g=\fso(n,\C)$ ($n\geq 6$). We find that $\La^x$  ($x\neq 0$)
is not a differential operator but instead is the left quotient of an order
$4$ algebraic differential operator  $D^x$ by an order $2$ invertible diagonalizable
operator. Precisely,
$\La^x=-\frac{1}{4}\frac{1}{E'(E'+1)}D^x$ where $E'$ is a positive shift of the
Euler vector field. So $\mu_x\star f$ is not local as an operator on $f$.
Thus $\star$ is not local.

The  differential operators  $D^x$  were constructed by us earlier (for this
purpose) in \cite{A-B:exotic}.
It would be very interesting now to find formulas for the operators $C_p(f,g)$ that
define $f\star g$. For $\g=\fsl(n+1,\C)$, 
some progress toward this is made in  \cite{me:RPn} 
using  results of  Lecomte and Ovsienko (\cite{LO}).  
Also we think that the method  of
Levasseur and Stafford (\cite{LSt}), which gave a new  
elegant construction of  our $D^x$ for $\g=\fsl(n+1,\C)$,  might  be  extended to
give the $\La^x$ and the $C_p(\cdot,\cdot)$.
These  approaches are based on  
the fact that $\RR$ identifies with the algebra of regular functions on $T^*(\CP^{n})$.

The star product  defines a representation $\pi$ of $\gog$  on $\RR$. We write this
out in Corollary \ref{cor:pi} using the $\La^x$.
In \S\ref{sec:9}, we show that $\star$ gives rise to  a
positive definite  hermitian inner product on $\RR$ compatible with $\pi$ and the
grading on $\RR$.  In this way, $\RR$ becomes the
Harish-Chandra module of a unitary representation of $G$  on the Hilbert space
completion $\cH=\widehat{\RR}=\widehat{\oplus}_{d=0}^\infty\RR^d$.
This quantizes $\Omin$, regarded as a real symplectic manifold, in the sense of
geometric quantization.  We compute the reproducing kernel of $\cH$ and deduce that
$\cH$ is a Hilbert space of holomorphic functions on $\Omin$.

It is a pleasure to thank Pierre Bieliavsky, Moshe Flato, Bert Kostant, Toby
Stafford, David Vogan and Ping Xu for helpful conversations. We also
warmly thank Brown University for their hospitality during the summers
of 1997 and 1998 when RKB was visiting there.
 
\section{Equivariant graded star  products on $\Sg/I$}  
\label{sec:2}

Let $\g$ be a complex semisimple Lie algebra.  The symmetric algebra 
$\cS=\Sg$  is the algebra of polynomial functions on $\g^*$. 
Then $\cS=\oplus_{d=0}^\infty\cS^d$ is a graded Poisson algebra in the
natural way, where $\{\cS^d,\cS^p\}\subseteq\cS^{d+p-1}$.
Let $I=\oplus_{d=0}^\infty I^d$ be a graded Poisson ideal in
$\cS$. We are most interested in the case when $I$ is the ideal $\I(\OO)$ of
functions vanishing on  a nilpotent coadjoint orbit $\OO$ in $\g^*$.
The term ``nilpotent" means that the corresponding adjoint orbit consists of
nilpotent elements; this happens if and only if $\OO$ is stable under dilations.

Let $\RR=\cS/I$ and  $\RR^d=\cS^d/I^d$. 
Then $\RR=\oplus_{d=0}^\infty\RR^d$ is again a graded Poisson algebra.
If $I=\I(\OO)$, then $\RR$ is the algebra of polynomial functions on the closure  $\Cl(\OO)$.
In the sense of algebraic geometry, $\Cl(\OO)$ is a closed complex algebraic subvariety of
$\g^*$ and $\RR$ is its algebra $R(\Cl(\OO))$ of regular functions. The elements $x\in\g$
define momentum functions $\mu_x$ in $\RR^1$ and  $\{\mu_x,\mu_y\}=\mu_{[x,y]}$.
The natural graded linear $G$-action on
$\RR$ corresponds to the $\g$-representation given by the operators $\Pmux$.

A \emph{graded star product} on $\RR$ is an associative  $\C[t]$-linear product
$\star$ on  $\RR[t]$ with the following properties. For  $f,g\in\RR$ we can write
$f\star g=\sum_{p=0}^\infty C_p(f,g)t^p$  and then
\begin{itemize}
\item[(i)]  $C_0(f,g)=fg $
\item[(ii)]  $C_1(f,g)-C_1(g,f)=\{f,g\}$
\item[(iii)]  $C_p(f,g)=(-1)^p C_p(g,f)$
\item[(iv)]  $C_p(f,g)\in\RR^{k+l-p}$  if    $f\in\RR^k$ and $g\in\RR^l$
\end{itemize}
Notice that (ii) and (iii) imply $C_1(f,g)=\half \{f,g\}$.
Axiom (iii) is called the \emph{parity} axiom.   

Given $\star$, we define a new noncommutative product on $\RR$ by
$f\circ g=f\star g\big|_{t=1}$.  Because of (iv), we can completely 
recover $\star$ from $\circ$.  
It is easy to see that (iii) amounts to the relation
$(f\circ g)^\al=g^\al\circ f^\al$ where $f\mapsto f^\al$ is the Poisson algebra
anti-involution of $\RR$ defined by $f^\al=(-1)^df$ if $f\in\RR^d$.

The star bracket is given by $[f,g]_\star=f\star g-g\star f$.
We say $\star$ is  \emph{$\g$-covariant} if  
$[\mu_x,\mu_y]_\star=t\mu_{[x,y]}$, and  
$\star$ is \emph{$G$-equivariant} (or 
strongly $\g$-invariant)    if  we have the much stronger relation
$[\mu_x, f]_\star=t\{\mu_x,f\}$.  
We say  that a   $G$-equivariant graded star product  on $\RR$ is  
an \emph{$G$-equivariant deformation quantization of $\RR$}.

Suppose $\star$ is a graded  $\g$-covariant star product  on $\RR$.
Let $\U=\Ug$ be the universal enveloping algebra of $\g$ equipped with its
canonical filtration $\{\U[d]\}_{d=0}^\infty$;   $\gr\U$ identifies naturally with
$\cS$. Then we have a noncommutative algebra homomorphism
$\Psi:\U\to\RR$ defined by 
$\Psi(x_1\cdots x_d)=\mu_{x_1}\circ\cdots\circ\mu_{x_d}$.
Then $\Psi$ is surjective in a filtered way, i.e., 
$\Psi(\U[p])= \oplus_{d=0}^p\RR^d$. The kernel of $\Psi$ is a $2$-sided ideal $J$ 
such that $\gr J=I$, and so  $\gr(\UgJ)$ identifies naturally with $\cS/I$.

Thus we get a vector space isomorphism  $\bq:\RR\to\UgJ$ defined by 
\begin{equation}\label{eq:bq} 
\bq(\mu_{x_1}\circ\cdots\circ\mu_{x_d})=x_1\cdots x_d+J
\end{equation}
Then  $\bq$ is a \emph{quantization map}, i.e.,
$\bq$ induces the identity maps $\RR^d\to\cS^d/I^d$.  We can recover $\circ$ 
from $\bq$  by the formula $f\circ g=\bq\inv((\bq f)(\bq g))$. Then
$\star$ is given by $f\star g=\bqt\inv((\bqt f)(\bqt g))$ where
$\bqt(f)=\bq(f)t^d$ if $f\in\RR^d$.

Let  $\tau$ be the algebra anti-involution of  $\U$ defined by $x^\tau=-x$; 
this is the so-called principal anti-automorphism.  The parity axiom (iii)  
implies that $J$ is stable under   $\tau$, so that $\tau$ descends to $\UgJ$,  and 
also $\bq(f^\al)=\bq(f)^\tau$. 
 
Clearly  $\star$ is $G$-equivariant  if and only if 
$\bq$ is $\g$-linear, i.e., $\bq(\{\mu_x,f\})=x\bq(f)-\bq(f)x$.
This amounts to $\bq$ being $G$-equivariant. In summary, this  discussion gives
\begin{prop}\label{prop:2} 
Suppose $\star$ is a graded  $G$-equivariant star product  on $\RR=\cS/I$.
Then we obtain in a canonical way a $2$-sided $\tau$-stable ideal $J$ in $\U$
and a $G$-equivariant quantization map $\bq:\RR\to\U/J$ given by 
\textup{(\ref{eq:bq})}.
\end{prop}

\section{Construction of $\star$ when $\OO=\Omin$}  
\label{sec:3}

From now on we assume that $\g$ is simple.    
Let $\Omin$ be the minimal non-zero nilpotent coadjoint orbit in $\g^*$.   So
$\Omin$  corresponds to the adjoint orbit of  highest root  vectors, or equivalently,
of highest weight vectors. We put $\RR=\cS/I$ where $I$ is the ideal of $\Omin$.

\begin{prop}\label{prop:e+u} 
Assume $\g$ is different from $\fsl(2,\C)$.   Then $\RR$ admits a  unique
$G$-equivariant graded star product $\star$.
\end{prop}

This strengthens the result in \cite{ABC} where they show that, if $\g$ is 
different from $\fsl(n+1,\C)$ for $n\ge 1$, then $\RR$ admits a $\g$-equivariant
graded star product which is unique up to equivalence of star products.
We need to exclude    $\g=\fsl(2,\C)$ because 
$\fsl(2,\C)$ admits infinitely many such star products 
(in natural bijection with $\C$).
\begin{proof} 
The discussion in \S\ref{sec:2} reverses easily to give a converse to 
Proposition \ref{prop:2}.   Precisely, if  $J$ is a $2$-sided
$\tau$-stable ideal of $\U$ such that $\gr J=I$ and $\bq:\RR\to\UgJ$  is a
$\g$-equivariant  quantization map such that  $\bq(f^\al)=\bq(f)^\tau$,
then the formula $f\star g=\bqt\inv((\bqt f)(\bqt g))$ defines a $\g$-equivariant
graded star product on $\RR$.
Thus it suffices to prove the following two statements.
\begin{itemize}
\item[(i)] There exists a unique $2$-sided ideal $J$ of $\U$ such that 
$\gr J=I$ and  $J^\tau=J$.
\item[(ii)]  For such $J$,
there exists a unique $G$-equivariant quantization map $\bq:\RR\to\UgJ$.
\end{itemize}
\vskip .25pc
\noindent
Notice that in (ii),  $\bq(f^\al)= \bq(f)^\tau$  follows automatically by   uniqueness.

To prove   (ii) we need only elementary facts about $\RR$  (see e.g., \cite{B-K:jams}).
The natural $G$-representation on $\RR$ is multiplicity free. 
One can get a very  quick abstract proof of (ii) just from this, but we will be 
more concrete.

$\RR^d$ is irreducible and carries  the $d$th Cartan power  $\g^{\boxt d}$
of the adjoint representation.   
Since the representation $\g^{\boxt d}$ occurs just
once in $\cS^d$, $I$ has  a unique graded $G$-stable complement 
$F=\oplus_{p=0}^\infty F^p$ in $\cS$; then $F$ identifies  with $\RR$.
We   define a vector space isomorphism 
$F\mapright{\bs}\,\U\mapright{}\,\UgJ$ where $\bs:\cS\to\U$ is the usual
symmetrization map; here we only assume that $\gr J=I$. 
Let $\bq$  be the corresponding  map from $\RR$ to $\UgJ$.   Then clearly
$\bq$ is a  $G$-equivariant quantization map.
If $\bh$ is another such map, then the composition $L=\bq\bh\inv$ satisfies:
$f\in\RR^d$ implies $L(f)=f+g$  where $g\in\RR^{d-1}$.
But also $L$ is $G$-linear  and so the $G$-decomposition of $\RR$ forces
$L(\RR^d)=\RR^d$. Thus $L(f)=f$.
 
The proof  of (i) breaks into two cases. If  $\g$ is 
different from $\fsl(n+1,\C)$, then as in \cite{ABC} we take $J$ to 
be the Joseph ideal constructed in \cite[\S5]{Jos}. 
We may characterize $J$ as the unique $2$-sided ideal in $\U$ whose associated
graded is $I$.  This is not the most familiar characterization, but it is immediate
from   the fact (\cite[Prop. 10.2]{Jos})  that $J$ is the
only completely prime   $2$-sided ideal such that  $\sqrt{\gr J}=I$, 
and the equality (\cite{Gar}) $\gr J=I$ .
Then uniqueness of $J$  implies that $J=J^\tau$.
 
Now suppose that $\g=\fsl(n+1,\C)$, $n\ge 1$.  
Let $\Dla(\CPn)$ be the algebra of global sections of the sheaf of twisted
differential operators acting on local sections of the $\la$th power of the 
canonical bundle on complex projective space; 
this makes sense for any complex number $\la$. We
have a natural algebra homomorphism $\Phi^{\la}:\U\to\Dla(\CPn)$. 
It is easy to write nice formulas for  the twisted vector fields $\Phi^\la_x$,
$x\in\g$,  in local coordinates on the big cell $\C^n$; see e.g., \cite{Tor}.

Let $\Jla$ be the kernel of $\Phi^\la$.  Then
$\Phi^\la$ is surjective and  $\gr\Jla=I$ (see \cite{BoBr}).
All $2$-sided  ideals $J$ in $\U$  with $\gr J=I$ arise in this way.
The principal anti-involution $\tau$ carries $\Jla$ to $J^{1-\la}$.
So $J=\Jh$ satisfies the two conditions in (i).
Assume $n\ge 2$. 
Then  we claim that $\Jla$ is $\tau$-stable iff $\la=\half$. 
To show this, we consider copies of the adjoint representation $\g$.

Since $\g$ appears (exactly) once in $\cS^2$, we see that $\g$ occurs twice in
$\U[2]$ and once in $\Jla_2=\Jla\cap\U[2]$.
The copy of $\g$ in $\cS^2$  corresponds, uniquely up to scaling, to a
$G$-linear map $r:\g\to\cS^2$,  $x\mapsto r^x$. 
Put $a^x=\bs(r^x)$. Then 
the copy of $\g$ in $\Jla_2$   consists  of elements $b^x=a^x+c_\la x$,
where $c_\la$ is some function of $\la$. A simple computation using the
formulas for $\Phi^\la_x$ mentioned above gives 
(for an appropriate choice of $r$)   $c_\la=\la-\half$. 
We have $\tau(a^x)=a^x$ while $\tau(x)=-x$.
So $b^x-\tau(b^x)=(2\la-1)x$. 
Thus, if $\la\neq\half$ then the unique copy of $\g$
in $\Jla_2$ is not  $\tau$-stable and consequently $\Jla$ is not  $\tau$-stable.
This proves the claim and finishes the proof of Proposition \ref{prop:e+u}.
\end{proof}
 
\begin{cor}\label{cor:sl2} 
Suppose $\g=\fsl(2,\C)$. Then   $\Jla=\Ker\Phi^\la$ 
\textup{(}$\la\in\C$\textup{)},  corresponds to a  $G$-equivariant graded star
product  $\starsub[\la]$ on $\RR$. All such star products arise in this way, and 
$\starsub[\la]=\starsub[\mu]$ iff $\mu=\pm(\la-\half)+\half$.
\end{cor}
\begin{proof} 
$\Jla$ is generated by a maximal ideal in the center of $\U$, 
it follows directly  that  $\Jla$ is $\tau$-stable and 
$\Jla=J^{\mu}$ iff $\mu=1-\la$.  
\end{proof}

\begin{prop}\label{prop:simple} 
In Proposition \textup{\ref{prop:e+u}},
the noncommutative algebra $\UgJ$ obtained by specializing $\star$ at $t=1$ is
a simple ring.
\end{prop}
\begin{proof} 
The Joseph ideal is maximal by \cite[Th.  7.4]{Jos}, and this means $\UgJ$ is
simple. If $\g=\fsl(n+1,\C)$, $n\ge 2$, then $\UgJ$ is isomorphic
to $\Dhalf(\CPn)$,  which is simple by   \cite{VdB}.
\end{proof}

\section{The operators $\La^x$}
\label{sec:4}

\begin{prop}\label{prop:4} 
The star product of a momentum function $\mu_x$, $x\in\g$, with an arbitrary
function  $f\in\RR$ is the three-term sum 
\begin{equation}\label{eq:mux*f} 
\mu_x\star f=\mu_xf+\thalf\Pmux[f]t+\La^x(f)t^2
\end{equation}
where  $\La^x$ are  operators on $\RR$. The $\La^x$    commute, are graded
of degree $-1$, and transform in the adjoint representation of $G$.
\end{prop}
\begin{proof} 
We have 
$\mu_x\star f=\mu_xf+\thalf\Pmux[f]t+\sum_{p=2}^\infty M^x_p(f)t^p$
where $M^x_p$ is graded of degree $-p$. Then
$x\otimes f\mapsto M^x_p(f)$ defines a $G$-linear map 
$M_p:\g\otimes\RR^d\to\RR^{d+1-p}$.
We know $\RR^d\simeq\g^{\boxt d}$  -- see the proof of Proposition
\ref{prop:e+u}(ii).  An easy fact about  representations (from highest weight
theory) is that if $\g^{\boxt k}$ appears $\g\otimes\g^{\boxt d}$ 
then $k$ lies in $\{d+1,d,d-1\}$. So $M_p=0$ if $p\ge 3$. Thus we get
(\ref{eq:mux*f}) where $\La^x=M_2^x$.

We have $(\mu_x\star f)\star\mu_y=\mu_x\star(f\star\mu_y)$.
Computing the coefficients of $t^4$, we find 
$\La^x\La^y(f)=\La^y\La^x(f)$. Computing the coefficients of $t^3$, we get
the relation
$[\eta^x,\La^y]=\La^{[x,y]}$ where $\eta^x=\Pmux$; so the $\La^x$
transform in the adjoint representation of $\g$.
\end{proof}

\begin{cor}\label{cor:La->star} 
\begin{itemize}
\item[(i)] The operators $\La^x$, $x\in\g$, completely determine $\star$.
\item[(ii)] The $\La^x$ generate a graded commutative
subalgebra $\A$ of  $\End\RR$ isomorphic to $\RR$.
\end{itemize}
\end{cor}
\begin{proof} 
(i) Once we know (\ref{eq:mux*f}), it is easy to compute
$\mu_{x_1}\cdots\mu_{x_k}\star f$ by induction on $k$.
(ii) This is easy, in fact $\La^{x_1}\cdots\La^{x_k}(f)$ is the coefficient of   
$t^{2d}$ in $\mu_{x_1}\cdots\mu_{x_k}\star f$.
Notice that $\A=\oplus_{d=0}^\infty\A^{-d}$ is graded in negative degrees, so that
$\A^{-d}$ corresponds to $\RR^d$.
\end{proof}

We have a representation $\pi$ of $\gog$ on $\RR$ defined by
$\pi^{x,y}(f)=\mu_x\circ  f-f\circ\mu_y$. 

\begin{cor}\label{cor:pi} 
The representation $\pi$ is irreducible and we have
\begin{equation}\label{eq:pixyphi=} 
\pi^{x,y}(f)=\mu_{x-y}f+\thalf\{\mu_{x+y},f\}+\La^{x-y}(f)
\end{equation}
\end{cor}
\begin{proof} 
$\pi$ is equivalent to the natural representation $\Pi$ of $\gog$ on $\UgJ$; indeed
$\bq$ is an intertwining map.   Proposition \ref{prop:simple} implies that $\Pi$ is
simple (and vice versa). 
\end{proof}

\begin{remark}\label{rem:kernels}  
Once we know the $\La^x$, we can construct $J$ directly as the kernel of the
algebra  homomorphism  $\U\to\End\RR$ defined by 
$x\mapsto\pi^{x,o}=\mu_x+\half\Pmux+\La^x$.
This is a noncommutative deformation of the fact that $I$ is the kernel of the 
algebra  homomorphism  $\cS\to\End\RR$ defined by $x\mapsto\mu_x$.
\end{remark}

The rest of this paper is devoted to computing the operators $\La^x$ when $\g$ is
classical. 

\section{The case $\g=\fsp(2n,\C)$}
\label{sec:5}
 Suppose  $\g=\fsp(2n,\C)$, $n\ge 1$. Let $\A$ be the Poisson algebra  
$\C[z_1,w_1,\dots,z_n,w_n]$ where  $\{z_i,z_j\}=\{w_i,w_j\}=0$ and
$\{z_i,w_j\}=\delta_{ij}$. 
We have a Poisson algebra grading $\A=\oplus_{k=0}^\infty\A^k$ where $\A^k$ is
the space of  homogeneous polynomials of total degree $k$.
Then $\A^2$ is a Lie subalgebra, and this is a model for  $\g$ (i.e., $\A^2$ is
isomorphic to $\g$). 
Moreover, $\A^{even}=\oplus_{k=0}^\infty\A^{2k}$ is a model for $\RR$.
The Moyal star product on $\A$ restricts to $\A^{even}$; in this way we get a
Moyal star product on $\RR$.

We find a strengthened version of  \cite[Prop.  6]{ABC}.
\begin{prop}\label{prop:sp} 
Let $\g=\fsp(2n,\C)$  $(n\ge 1)$.  
The Moyal star product   is a $G$-equivariant graded star product on $\RR$.
If   $n\ge2$, then it corresponds to the Joseph ideal;
if $n=1$, then it  corresponds to the ideal $J^{\frac{1}{4}}$. 

The  $\La^x$ are   order  $2$ algebraic differential  operators and 
\begin{equation}\label{eq:sp} 
\La^{z_iz_j}=\frac{1}{4}\frac{\del^2}{\del w_i\del w_j},
\qquad
\La^{w_iw_j}=\frac{1}{4}\frac{\del^2}{\del z_i\del z_j},
\qquad
\La^{z_iw_j}=-\frac{1}{4}\frac{\del^2}{\del w_i\del z_j}
\end{equation} 
\end{prop}

\section{Computation of $\La^x$}
\label{sec:6}
We assume from now on that   $\g$ is a classical complex simple 
Lie algebra different from $\fsp(2n,\C)$, $n\geq 1$. This falls into  two cases:
(I)  $\g=\fsl(n+1,\C)$ where  $n\geq 2$, or
(II) $\g=\fso(n,\C)$ where  $n\geq 6$.
It turns out that we can deal with both cases simultaneously by simply 
by  introducing a parameter $\ep$ and setting $\ep=0$ in  (I) 
or $\ep=1$  in (II).  We set $\pep=p+\ep$ and $\pmep=p-\ep$.

We put $G=SL_n(\C)$ in (I) or $G=Spin_n(\C)$ in (II).
Notice  that there is one coincidence between  (I) and (II), namely
$\g=\fsl(4,\C)=\fso(6,\C)$. 

We define  $m$ by  $\dim\,\Omin=2m+2$;
so $m=n-2$ in (I) or $m=n-4$ in (II).
Let $X,h,Y$ be a triple in $\g$ such that $X\in\Omin$ and 
$[X,Y]=h$, $[h,X]=2X$,   $[h,Y]=-2Y$.
Then $h$ is semisimple and  $Y\in\Omin$. 
In this same setting  we proved    
\begin{thm}[\cite{A-B:exotic}]  \label{thm:exotic}  
Let $\D_{4;-1}(\Omin)$ denote the space of  algebraic differential operators $D$
on $\Omin$ such that $D$ has order at most $4$ and $D$ is  graded of degree $-1$,
i.e., $D(\RR^p)\subseteq\RR^{p-1}$.  

Then $\D_{4;-1}(\Omin)$ contains a unique copy of the adjoint representation of
$G$. In other words, there is a non-zero
$G$-equivariant complex linear map 
$\g\to\D_{4;-1}(\Omin)$, $x\mapsto D^x$, and this map is
unique up to scaling.
For $x\neq 0$,  $D^x$ has order exactly $4$.
  
We can normalize the map  $x\mapsto D^x$  so that, for $p\ge 0$, 
\begin{equation}\label{eq:DYmuxp} 
D^Y(\mu_X^p)=\ga_p\mu_X^{p-1}
\end{equation}
where  $\ga_p=p(p+\tfrac{m-1}{2})\pep(\pmep+\tfrac{m}{2})$.
For $p\ge 1$,  $D^Y(\mu_X^p)\neq 0$.

Finally,   the map $x\mapsto D^x$ extends naturally to an
the operators $D^x$ generate a graded commutative subalgebra of $\D(\Omin)$
which is isomorphic to $\RR$.  Thus for $f\in\RR$ we have the operator $D^f$,
where  $D^fD^g=D^{fg}$ and $D^{\mu_{x}}=D^x$.
\end{thm}
\begin{proof} 
This is a summary of the following results in  \cite{A-B:exotic}:  
Theorems 3.2.1 and 3.2.3, Corollary 3.2.4,  Propositions 4.2.3 and  4.3.3, and Corollary 3.2.5.  
\end{proof}


\begin{remarks}\label{rems:6} 
(i) If  $\g=\fsl(4,\C)=\fso(6,\C)$, then we can equally well choose 
$\ep=0$ or $\ep=1$ in computing $\ga_p$.
We of course end up with the same final answer.

\noindent (ii) $D^x$ defines an algebraic differential operator on $\Cl(\Omin)$.
\end{remarks}

Let $E$ be the  Euler vector field on $\Omin$ so that $E$ operates on $\RR$ and
$\RR^d$ is its $d$-eigenspace.   
We put $E'=E+\tfrac{m+1}{2}$. Notice that
$E'$ is diagonalizable  on $\RR$ with
positive spectrum, and so $E'+k$ is invertible for any $k\ge 0$. 

\begin{thm}\label{thm:Lax=} 
For $x\in\g$ we have
$\La^x= -\,\dfrac{1}{4E'(E'+1)}D^x$.
\end{thm}
\begin{proof} 
This occupies \S\ref{sec:7}.
\end{proof}

We found  this formula for $\La^x$  because we expected this shape
$\La^x= P\inv D^x$ where $P$ is a quantization of $4\la^2$ and $\la$ is the
symbol of $E$; see \cite[\S1]{A-B:exotic}.

\begin{remark}\label{rem:sp}  
We can fit the case $\g=\fsp(2n,\C)$ discussed in \S\ref{sec:5} into this framework
formally  by putting   
$D^x=-4E'(E'+1)\La^x$ where the $\La^x$ were given in (\ref{eq:sp})
and again $E'=E+\tfrac{m+1}{2}$ for $m=\half\dim\,\Omin-1=n-1$.
Then the formula   $ D^Y(\mu_X^p)=\ga_p\mu_X^{p-1}$  in Theorem \ref{thm:exotic}
still holds if we compute $\ga_p$ for   $\ep=-\half$. 
Here we may choose $X=-\half w_1^2$, $Y=\half z_1^2$, $h=z_1w_1$.
\end{remark}

\section{Proof of Theorem \ref{thm:Lax=}}
\label{sec:7} 

\begin{lem}\label{lem:phiD} 
We have $\La^x=\phi D^x$  where $\phi$ is a linear operator on $\RR$ given by
scalars $\phi_d$, $d\ge 0$, so that  $\phi(f)=\phi_{d-1}f$ if $f\in\RR^d$.
The scalars $\phi_d$ are   unique.
\end{lem}
\begin{proof} 
Let $p\ge 1$.
We have two $G$-linear maps  $\g\otimes\RR^p\to\RR^{p-1}$ defined by
$\al_p(x\otimes f)=\La^x(f)$ and $\be_p(x\otimes f)=D^x(f)$. 
These must be proportional because   
$\Hom_G(\g\otimes\g^{\boxt p},\g^{\boxt(p-1)})$ is $1$-dimensional. 
We know that $\be_p$ is non-zero by  Theorem \ref{thm:exotic}.
So  there is a unique scalar $\phi_{p-1}$ such that $\al_p=\phi_{p-1}\be_p$.
\end{proof}
 
At this point,  there is no guarantee that $\phi_p$ will be a nice function of
$p$, in the sense that $\phi$ is a reasonable function of  $E$.  
But Theorem \ref{thm:Lax=} asserts  $\phi=-\frac{1}{4}\frac{1}{E'(E'+1)}$.
 
To prove this, 
we will write down a series of recursion relations for the $\phi_p$.
To derive the recursions, we start with the bracket relation
$[\pi^{x,-x},\pi^{y,-y}]=\pi^{z,z}$ where $z=[x,y]$.  
By (\ref{eq:pixyphi=})   we have
$\pi^{x,-x}=2\mu_x+2\La^x$ and $\pi^{z,z}=\eta^z$
where $\eta^z=\{\mu_z,\cdot\}$.  Since  the operators $\La^x$, like the  
$\mu_x$, commute among themselves, we get
\begin{equation}\label{eq:bracket} 
[\mu_x,\La^y]+[\La^x,\mu_y]=\tfrac{1}{4}\eta^{[x,y]}
\end{equation}
We  choose  $x=X$ and $y=Y$  so that $[x,y]=h$.  
Writing  $\La^x=\phi D^x$ and applying the  operator identity
(\ref{eq:bracket}) to a test function $f\in\RR^p$, $p\ge 1$, we find
\begin{equation}\label{eq:key}
\phi_{p-1}\mu_XD^Y(f)-\phi_pD^Y(\mu_Xf)+\phi_pD^X(\mu_Yf)-
\phi_{p-1}\mu_YD^X(f)=\tfrac{1}{4}\eta^h(f).
\end{equation}
The recursions will arise by evaluating this for   $f=\mu^s_X\mu^t_Y$,
with  $s+t=p$.
 
Before we can  write down the  recursions, we  need some auxiliary computations,
provided by the next result.
(Unfortunately,  (\ref{eq:DYmuxp}) is not sufficient to  determine all $\phi_p$.) 

\begin{lem}\label{lem:albe}
For $s,t\geq 0$ we have
\begin{eqnarray}
D^Y(\mu^s_X\mu^t_Y)=\al_{s,t}\mu_X^{s-1}\mu^t_Y+
\be_{s,t}\mu_X^{s-2}\mu_Y^{t-1}\mu_h^2\label{eq:DYprod}\\[6pt] 
D^X(\mu^t_X\mu^s_Y)=
\alpha_{s,t}\mu^t_X\mu^{s-1}_Y+\beta_{s,t}\mu^{t-1}_X
\mu^{s-2}_Y\mu^2_h\label{eq:DXprod}
\end{eqnarray}
where  
$\alpha_{s,t}=\gamma_s+\half st(2s+t+m)$ and 
$\beta_{s,t}=-\frac{1}{4}(s-1)st(2s+t+m)$.
 \end{lem}
\begin{proof} 
We have to go back into our explicit construction of $D^Y$ in
\cite[\S4]{A-B:exotic}. We worked over the Zariski open dense set
$\Oreg=(\mu_Y\neq 0)$   in $\Omin$. We constructed  $D^Y$  as the quotient
$D^Y=\frac{1}{\mu_Y}S$ where  $S$ is a certain differential operator on $\Oreg$.
More precisely, $S=\frac{1}{4}(T-q(\eta^Y)^2)$   where 
$q=(E+\frac{m}{2}+\ep)(E+\frac{m}{2}-\ep)$  and 
$T$ is an explicit noncommutative
polynomial in some vector fields on $\Oreg$ which annihilate 
$\mu_Y$. Also $\eta^Y$ annihilates $\mu_Y$. It follows that for
any $g\in R(\Oreg)$ we have
$T(g\mu^t_Y)=T(g)\mu^t_Y$ and  $\eta^Y(g\mu^t_Y)=\eta^Y(g)\mu^t_Y$.
 
Now we can  compute $D^Y(\mu^s_X\mu^t_Y)$. We have
$D^Y=A-B$ where  
$A=\frac{1}{4\mu_Y}T$ and $B=\frac{1}{4\mu_Y}q(\eta^Y)^2$.
Then  we find
$D^Y(g\mu^t_Y)=D^Y(g)\mu^t_Y+B(g)\mu^t_Y-B(g\mu^t_Y)$.  
Let  $g=\mu^s_X$. Then 
$D^Y(\mu^s_X)=\gamma_s\mu^{s-1}_X$ by (\ref{eq:DYmuxp}). 
Also, since $\eta^Y(\mu_X)=\mu_{[Y,X]}=-\mu_h$ and
$\eta^Y$ is a  vector field we find (as in \cite[(67)]{A-B:exotic})
\begin{equation}\nonumber 
(\eta^Y)^2(\mu^s_X)=\left(-2s\mu_X\mu_Y+s(s-1)\mu_h^2
\right)\mu^{s-2}_X
\end{equation} 
Using this we find 
\begin{equation}\nonumber
B(\mu^s_X\mu^t_Y)=
\tfrac{1}{4}sq_{s+t}\left(-2\mu_X+
(s-1)\mu^{-1}_Y\mu_h^2\right)\mu^{s-2}_X\mu^{t}_Y
\end{equation}
where  $q_p=(p+\frac{m}{2}+\ep)(p+\frac{m}{2}-\ep)$.  Now we   obtain
\begin{equation}\nonumber  
\begin{array}{ll}
D^Y(\mu^s_X\mu^t_Y) 
&=\gamma_s \mu^{s-1}_X\mu^t_Y-\frac{1}{4}
s(q_{s+t}-q_s)\left(-2\mu_X+(s-1)\mu^{-1}_Y\mu_h^2\right)
\mu^{s-2}_X\mu^{t}_Y\\[6pt]
&=\alpha_{s,t}\mu^{s-1}_X\mu^t_Y+\beta_{s,t}
\mu^{s-2}_X\mu^{t-1}_Y\mu^2_h
\end{array}
\end{equation} 
where  $\alpha_{s,t}=\gamma_s+\half s (q_{s+t}-q_s)$ and 
$\beta_{s,t}=-\frac{1}{4}(s-1)s(q_{s+t}-q_s)$.
This proves   (\ref{eq:DYprod}).

We can prove (\ref{eq:DXprod}) by applying  a certain automorphism. Let
$\chi:SL(2,\C)\to G$ be the Lie group homomorphism
corresponding to the Lie algebra inclusion $\fs\to\g$ where $\fs$
is the span of $X$, $h$, and $Y$. The adjoint action of  
$\chi\left(\begin{smallmatrix}\phantom{-}0& 1\\ -1 &0\end{smallmatrix}\right)$
defines a Lie algebra automorphism  $\vth$ of $\g$. Then
$\vth(X)=-Y$, $\vth(Y)=-X$ and $\vth(h)=-h$. Clearly 
$\vth$ preserves $\Omin$ and hence induces algebra
automorphisms of $\RR$ and of $\D(\Cl(\Omin))$ which we again
call $\vth$.  Then  $\vth(\mu_x)=\mu_{\vth(x)}$
and $\vth(D^x)=D^{\vth(x)}$.  Now applying $\vth$ to 
(\ref{eq:DYprod}) we get (\ref{eq:DXprod}).
\end{proof}
\begin{remark}\label{rem:use_LS}  
For  $\g=\fsl(n+1,\C)$,  these calculations  become much easier   
if  we use the  formulas  for $D^x$  found in \cite{LSt}. 
But there are no  such formulas known  when $\g=\fso(n,\C)$.
\end{remark}

Now we can obtain  the recursions by
plugging $f=\mu^s_X\mu^t_Y$, where $p=s+t$, into (\ref{eq:key}). We evaluate
using  (\ref{eq:DYprod}), (\ref{eq:DXprod})  and the fact $\eta^h(f)=2(s-t)f$. 
The result only involves two functions, namely $f$  and 
$g=\mu^{s-1}_X\mu^{t-1}_Y\mu^2_h$.  We find, for  $s,t\ge 0$,
\begin{equation}\nonumber%
\begin{array}{l}
\phi_{p-1}\left[(\al_{s,t}-\al_{t,s})f+(\be_{s,t}-\be_{t,s})g\right]  
-\phi_p\left[(\al_{s+1,t}-\al_{t+1,s})f+(\be_{s+1,t}-\be_{t+1,s})g\right]  \\[6pt]
=\half(s-t)f 
\end{array}
\end{equation}
Equating coefficients of $f$ and  $g$  we obtain the two recursions
\begin{eqnarray}
\phi_{p-1}(\alpha_{s,t}-\alpha_{t,s})-\phi_p(\alpha_{s+1,t}
-\alpha_{t+1,s})&=&\textstyle\half(s-t),
\label{eqR1}\\[5pt]
\phi_{p-1}(\beta_{s,t}-\beta_{t,s})-\phi_p(\beta_{s+1,t}-
\beta_{t+1,s})&=&0\label{eqR2}
\end{eqnarray}
Both recursions are valid for $s,t\ge 1$, since $f$ and $g$ are linearly independent
functions on  $\Omin$.  Moreover   (\ref{eqR1}) is valid for all
$s,t\ge 0$,    since $\be_{i,j}=0$ if $i=0$, $i=1$ or  $j=0$.
 
First we  consider   (\ref{eqR2}). 
Our formula for $\be_{s,t}$ in Lemma \ref{lem:albe} yields 
\begin{equation}\label{eq:bes} 
\begin{array}{ccc} 
\beta_{s,t}-\beta_{t,s}&=&-\tfrac{1}{4}st(s-t)(2s+2t+m-1)
\\[4pt]
\beta_{s+1,t}-\beta_{t+1,s}&=&-\tfrac{1}{4}st(s-t)(2s+2t+m+3)
\end{array}
\end{equation}
For $p\geq 3$   we can write $p=s+t$ with $s,t\ge 1$ and
$s\neq t$. Then (\ref{eqR2})  and  (\ref{eq:bes}) give  
\begin{equation}
\label{eqR3}
\phi_{p}=\frac{2p+m-1}{2p+m+3}\phi_{p-1},\qquad p\ge 3
\end{equation}
This is  a very simple recursion  with solution 
\begin{equation}\label{eqR4}
\phi_p=\phi_2\frac{(m+5)(m+7)}{(2p+m+1)(2p+m+3)},
\qquad p\ge 2
\end{equation}

Our  aim is to prove $\phi=-\frac{1}{4}\frac{1}{E'(E'+1)}$, which amounts to  
$\phi_p=-\frac{1}{4d_p(d_p+1)}$, $p\ge 0$, where $d_p=p+\frac{m+1}{2}$.
So we are pleased that (\ref{eqR4}) gives 
\begin{equation}\label{eqR7}
\phi_p=\frac{\om}{4d_p(d_p+1)},\qquad p\ge 2
\end{equation}
where $\om$ is the constant  $(m+5)(m+7)\phi_2$.

To determine  $\phi_p$ at
$p=0,1,2$, we   implement  (\ref{eqR1}) for $t=0$ and $p=s$. Since
 $\alpha_{p,0}=\ga_p$ and $\al_{0,p}=0$  we get 
\begin{equation}\label{eqR5}
\phi_{p-1}\ga_p-\phi_p(\ga_{p+1}-\al_{1,p})=
\thalf p,\qquad p\ge 1
\end{equation}
To use this, we observe $\ga_p=pd_p\nu_p$, $p\ge 0$,
where   $\nu_p=\pep(\pmep+\tfrac{m}{2})$.
After a little  work,  we find 
$\ga_{p+1}-\al_{1,p}=p(d_p+1)(\nu_p+2d_p)$.   Now (\ref{eqR5}) gives 
\begin{equation}\label{eqR6}
\phi_p=\frac{\nu_p(d_p-1)\phi_{p-1}-\thalf}{(\nu_p+2d_p)(d_p+1)},
\qquad  p\geq 1
\end{equation} 
If we put  $\la_p=4d_p(d_p+1)\phi_p$ ($p\ge 0$)  this simplifies nicely to give
\begin{equation}\label{eqR9} 
\nu_p(\la_p-\la_{p-1})=-2d_p(\la_p+1),\qquad   p\ge 1 
\end{equation}
Plugging in $\la_p=\om$ for $p\ge 2$, we get
$\om=-1$, and then $\la_1=\la_0=-1$. 
Thus, for all $p\ge 0$, $\la_p=-1$ and so  $\phi_p=-\frac{1}{4d_p(d_p+1)}$.

\begin{remark}\label{rem:new}  
In this proof, we only used the fact that  there exists
some  $J$ such that $\gr J=I$ and $J=J^\tau$. But now (see Remark
\ref{rem:kernels}) we can recover $J$
as the kernel of the algebra  homomorphism  $\U\to\End\RR$ defined by 
$x\mapsto\mu_x+\half\Pmux-\frac{1}{4E'(E'+1)}D^x$.
This gives a different proof that $J$ is unique.
\end{remark}  

\section{Consequences of Theorem \ref{thm:Lax=}}
\label{sec:8}

We may rescale the complex Killing form
$\Kill{\cdot}{\cdot}$ of $\g$ so  that  $\Kill{X}{Y}=\half$.

\begin{cor}\label{cor:Lax} 
\begin{itemize}
\item[(i)]  We have $\La^Y(\mu_X^p)=\ze_p\,\mu_X^{p-1}$ where
$\ze_p=-\frac{\ga_p}{(2p+ m+1)(2p+ m+3)}$
\item[(ii)] The map $\La^x:\RR^p\to\RR^{p-1}$ is non-zero if
$p\ge 1$ and $x\neq 0$.
\item[(iii)] $\La^x(y)=c\Kill{x}{y}$  where $c$ is a non-zero scalar; in fact
$c=2\ze_1$.
\end{itemize}
\end{cor}
\begin{proof} 
(i)  is  immediate from (\ref{eq:DYmuxp}).
This gives (ii) if $x=Y$ (since $\ga_p\neq 0$ if $p\ge 1$).
Since  the $\La^x$ transform in the adjoint representation,  
we  get (ii) for all $x$.
Finally (iii) follows because the map 
$\g\otimes\g\to\C$, $x\otimes y\mapsto\La^x(\mu_y)$, is $G$-invariant and so
must be a multiple $c\Kill{\cdot}{\cdot}$ of our normalized  Killing form
(see  \S\ref{sec:7}). Then $c$ is non-zero by (ii);  choosing $x=Y$
and $y=X$ we find $c=2\ze_1$.
\end{proof}

\begin{cor}\label{cor:fails} 
For $x\neq 0$, $\La^x$ fails to be a differential operator on $\Omin$.
In fact, neither factor $E'$ nor $E'+1$ left divides $D^x$.
\end{cor}
\begin{proof} 
Suppose one of  $E'$ or $E'+1$ left divides $D^x$ so that  the quotient
is a differential operator  $A^x$  on $\Omin$.  Since $D^x$ has order
$4$ (Theorem \ref{thm:exotic}),    $A^x$ has order  $3$.
But then the $A^x$ span a copy of the adjoint representation in 
$\D_{4;-1}(\Omin)$ which is different from the copy spanned by the $D^x$.
This  contradicts   uniqueness in Theorem \ref{thm:exotic}.
\end{proof}    

Notice that the corollary implies that $\La^x$ fails to be a differential operator on
$\RR$ (since otherwise $\La^x$ would be a differential operator on $\Cl(\Omin)$).

\begin{remark}\label{rem:psd} 
Theorem \ref{thm:Lax=}   suggests that 
$C_2(\mu_x,\cdot)=\La^x$ is ``pseudo-differential" in some sense.
This is  different in character from the often cited  example of
``pseudo-differential" star product found  \cite[\S9, page 124]{Fron} for
coadjoint orbits of the Euclidean group $E(2)$.  There Fronsdal obtains a star
product    where the operator  $f\mapsto\mu_x\star f$ is   an infinite series
of differential operators $C_k(\mu_x,\cdot)t^k$ with increasing order.
\end{remark}

There is a unique (up to scaling) casimir in $\U[2]$, namely
$Q=\sum_{i=1}^N{x_i}^2$,  where  $\{x_i\}_{i=1}^N$ is a basis of $\g$ such that  
$\Kill{x_i}{x_j}=\delta_{ij}$.  We next compute how $Q$ acts on $\UgJ$
with respect to the left multiplication action of $\U$.
 
\begin{cor}\label{cor:} 
$Q$ acts on $\UgJ$ by the scalar  
$s=-\dfrac{(1+\ep)(m+2-2\ep)}{4(m+3)}\dim\,\g$.
\end{cor}
\begin{proof} 
The function $\sum_{i=1}^N\mu_{x_i}^2$ is $G$-invariant and so vanishes on
$\Omin$.   Now (\ref{eq:mux*f}) and Corollary \ref{cor:Lax} give
$\sum_{i=1}^N\mu_{x_i}\star\mu_{x_i}=
\sum_{i=1}^N(\mu_{x_i}^2+c\Kill{x_i}{x_i})=2\ze_1N$.
This means (see \S\ref{sec:2}-\ref{sec:3}) that $Q-2\ze_1N$ lies in $J$, and so $Q$
acts by   $2\ze_1N$.  
\end{proof}

\begin{remark}\label{rem:conj}  
We conjecture that for the $5$ exceptional simple Lie
algebras, $\La^x$ again has the form $-\frac{1}{4E'(E'+1)}D^x$  where $D^x$ are
some (as yet unknown) order $4$ algebraic differential operators on $\Omin$.
\end{remark}
 
\section{Hermitian inner  product on $\RR$}
\label{sec:9}
We assume that $\g$ is a complex simple Lie algebra different from $\fsl(2,\C)$.
Let $U$ be a maximal compact subgroup of $G$. Let $\sig$ be the  corresponding
Cartan involution  of $\g$; so $\sig$ is $\C$-antilinear.  
Then $\gs=\{(x,\sig(x))\,|\, x\in\g\}$ is a  real form of $\gog$.
We have a $U$-invariant $\C$-antilinear algebra involution 
$f\mapsto\of$ on $\RR$ defined by $\of(z)=\ovl{f(\sig(z))}$; see \cite[\S2.3]{A-B:exotic}.

\begin{thm}\label{thm:herm} 
The formula 
\begin{equation}\label{eq:f|g} 
(f|g)=\mbox{constant term in } f\circ\og  
\end{equation} 
defines a $U$-invariant positive definite  hermitian inner product on $\RR$, with $(1|1)=1$.
In addition  $(\cdot|\cdot)$ is  $\gs$-invariant, i.e.,  the operators $\pi^{x,\sig(x)}$,
$(x\in\g)$, are skew-adjoint.
\end{thm}
\begin{proof} 
The pairing  $(\cdot|\cdot)$ is clearly sesquilinear and $U$-invariant   with $(1|1)=1$.
It follows by $U$-invariance that  $\RR^j$ is orthogonal to $\RR^k$ if $j\neq k$.
Now to show $(\cdot|\cdot)$ is hermitian positive definite, it suffices to show that each number
$\ns{\mu_Y^p}=(\mu_Y^p|\mu_Y^p)$, $p\ge 1$,  is   positive. 
We will use Remark \ref{rem:sp} so that  we can treat all cases simultaneously.
We may assume now, by rechoosing $(X,h,Y)$ if
needed, that $\sig(Y)=-X$; see   \cite[\S2.3]{A-B:exotic}. Then by Corollary \ref{cor:pi}, Theorem 
\ref{thm:Lax=} and   (\ref{eq:DYmuxp})   we  have
\begin{equation}\label{eq:long} 
\ns{\mu_Y^p}=(-1)^p(\La^Y)^p(\mu_X^p)=(\tfrac{1}{4E'(E'+1)}D^Y)^p(\mu_X^p)
=\textstyle\prod_{i=1}^p\frac{\ga_i}{(2i+ m-1)(2i+ m+1)}
\end{equation}
This number is positive, since  $\ga_p$ is positive by  Theorem \ref{thm:exotic}.

For $f\in\RR$  we have, by Corollary \ref{cor:La->star}, the operator $\La^f$ on $\RR$,
where  $\La^f\La^g=\La^{fg}$ and $\La^{\mu_{x}}=\La^x$. Then plainly
\begin{equation}\label{eq:short} 
(f|g)= \mbox{constant term in } \La^f(\og)
\end{equation}
This formula easily implies that the adjoint of (ordinary) left multiplication by $\mu^x$ is
$\La^{\sig(x)}$. Hence the operators  
$\mu_{x} -\La^{\sig(x)}$ are skew-adjoint. But also the operators 
$\{\mu_{x+\sig(x)},\cdot\}$ are skew-adjoint since they correspond to the action of $U$.
Thus, using   (\ref{eq:pixyphi=}), we see the  operators
$\pi^{x,\sig(x)}$  are all skew-adjoint.
\end{proof}
Notice that, since $\pi$ is irreducible by Corollary \ref{cor:pi},  $(\cdot|\cdot)$ is the unique
$\gs$-invariant hermitian pairing on $\RR$ such that $(1|1)=1$.

\begin{cor}\label{cor:unit} 
The operators $\pi^{x,\sig(x)}$ on $\RR$ exponentiate to give a unitary representation of $G$
on the Hilbert space direct sum $\cH=\widehat{\oplus}_{d=0}^\infty\RR^d$.
Then $\RR$  is the Harish-Chandra module of this unitary representation.
\end{cor}
\begin{proof} 
This follows by a theorem of Harish-Chandra since 
$\RR$ is an admissible  $(\gog,G)$-module where
$\gog$ acts by $\pi$ and $G$ acts  corresponding to the operators $\{\mu^x,\cdot\}$.
\end{proof}
This quantizes $\Omin$ in the sense of geometric
quantization and the orbit method. We note that the shift from $E$ to $E'$ can be explained
by half-forms  in the same way as in \cite[Prop. 5]{me:gqmin}.

\begin{cor}\label{cor:repro} 
The unitary representation of $G$ on $\cH$ admits a reproducing kernel $\cK$.
Explicitly, $\cK$ is the
function $\cK(x,y)$ on  $\Omin\times\Omin$ given by the hypergeometric function
\begin{equation}\label{eq:cK=} 
\cK= \phantom{x}_1F_2\left(\frac{m+3}{2}; 1+\ep, 1-\ep+\frac{m}{2}; 2T\right)
\end{equation}
where   $T(x,y)=-\Kill{x}{\sigma(y)}$.  So $\cK(x,y)$ is holomorphic in $x$ and
anti-holomorphic in $y$.
Consequently, $\cH$ is a Hilbert space of
holomorphic functions on $\Omin$.
\end{cor}

\begin{proof} 
Going back to (\ref{eq:long}), we find
\begin{equation}\label{eq:new}
\ns{\mu_Y^p}=\frac{p!(1+\ep)_p (1-\ep+\tfrac{m}{2})_p }{4^p(\frac{m+3}{2})_p}
\end{equation}
where we are using the classical notation $(a)_p=a(a+1)\cdots(a+p-1)$.
By definition, 
$\cK=\sum_{i=0}^\infty f_i\otimes\of_i$ where $f_0, f_1,\dots$ is an orthonormal basis
of $\RR$ with respect to $(\cdot|\cdot)$. 
On the other hand, $T= \sum_{i=0}^N s_i\otimes\ovl{s}_i$ where
$s_0,\dots,s_N$ is an orthonormal basis of $\RR^1$ with respect to the hermitian inner
product $\langle\mu_x|\mu_y\rangle=-\Kill{x}{\sigma(y)}$.
This is positive definite  since $\langle\mu_Y|\mu_Y\rangle=\Kill{Y}{X}=\half$.
It follows, as in \cite[\S8]{me:gqmin}, that   
\begin{equation}\label{eq:cK==} \nonumber
\cK=\sum_{p=0}^\infty \,\frac{1}{\ns{\mu_Y^p} }\,\left(\frac{T}{2}\right)^p 
\end{equation}
So   (\ref{eq:new})  gives (\ref{eq:cK=}).  
\end{proof}

\begin{remark}\label{rem:last}  
In the case $\g=\fsp(2n,\C)$,  $\cH$ is just the   classical Fock space 
of even holomorphic functions $f(z_1,w_1,\dots,z_n,w_n)$ with reproducing kernel
$\cK=\cosh(2\psi)$  where $\psi=\sum_{i=1}^n (|z_i|^2+|w_i|^2)$.
Indeed, $T=\half\psi^2$ and the hypergeometric series  collapses  to 
\begin{equation}\label{eq:collapse} 
\phantom{.}_1F_2\left(\frac{m+3}{2}; \half, \frac{m+3}{2}; 2T\right)
=\cosh(\sqrt{8T})=\cosh(2\psi) 
\end{equation}  
\end{remark}

\bibliographystyle{amsalpha}

\end{document}